\documentclass[letterpaper, reqno,11pt]{article}
\usepackage{amsfonts}
\usepackage{bbm}
\usepackage{amsmath,amssymb}
\usepackage{mathrsfs}
\usepackage{hyperref}
\usepackage{graphicx}
\usepackage{float}
\usepackage{amsmath,amscd,amsbsy,amssymb,latexsym,url,bm,amsthm}
\usepackage{mathrsfs}
\usepackage{hyperref}
\usepackage{graphicx}
\usepackage{float}

\usepackage{amsmath,amscd,amsbsy,amssymb,latexsym,url,bm,amsthm}
\usepackage[vlined,boxed,commentsnumbered,ruled]{algorithm2e}
\usepackage{epsfig,graphicx,subfigure}
\usepackage{enumitem,balance,mathtools}
\usepackage{wrapfig}
\usepackage{mathrsfs, euscript}
\usepackage[usenames]{xcolor}
\usepackage{caption}
\usepackage{lipsum}
\usepackage{siunitx}		%数字单位科学计数法等
\usepackage{bibentry}
\usepackage{tikz} 		%画图

\newtheorem{thm}{Theorem}[section]
\newtheorem{coro}[thm]{Corollary}
\newtheorem{lemma}[thm]{Lemma}
\newtheorem{prop}[thm]{Proposition}

\newtheorem{rem}[thm]{Remark}

\numberwithin{equation}{section}\allowdisplaybreaks

\newcommand{\besov}{B^{s}_{p,q}}
\newcommand{\besovo}{B^{s}_{p_{0},q_{0}}}
\newcommand{\besovon}{B^{s}_{p_{1},q_{1}}}

\newcommand{\triebel}{F^{s}_{p,q}}

  
\makeatletter

\def\leq{\leqslant}

\def\eps{\varepsilon}
\def\leq{\leqslant}
\def\geq{\geqslant}
\def\real{{\mathbb{R}}}			%实数
\def\FF{{\mathscr{F}^{-1}}}		%Fourier 逆变换
\def\F{{\mathscr{F}}}			%Fourier 变换
\def\Z{{\mathbb{Z}}}			%整数
\def\sch{\mathscr{S}}

\def\mpq{M_{p,q}}

\usepackage[percent]{overpic}
\usepackage[margin=1.0in]{geometry}
\usepackage{color,latexsym,amsmath,amssymb, amsthm, graphicx, enumitem}

\newcommand{\norm}[1]{\left\|#1\right\|}		%模竖线
\newcommand{\abs}[1]{\left|#1\right|}			%绝对值
\newcommand{\set}[1]{\left\{#1\right\}}			%集合括号
\newcommand{\inner}[1]{\left\langle#1\right\rangle} %内积括号
\newcommand{\rmnum}[1]{\uppercase\expandafter{\romannumeral#1}}  %罗马数字

\begin{document}
	\pagenumbering{arabic}	
	\title{Sharp embedding between Besov-Triebel-Sobolev spaces and modulation spaces }
	\author{Yufeng Lu\thanks{School of Mathematics Sciences, Peking University, Beijing 100871, PR China, luyufeng@pku.edu.cn}}
	\date{}
	\maketitle
%	\tableofcontents	%制作目录
%	\newpage

	\begin{abstract}
    The embedding relations between Besov-Triebel-Sobolev spaces and modulation spaces are determined explicitly. We extend the results in \cite{Sugimoto2007dilation},\cite{Wang2007Frequency},\cite{Kobayashi2011inclusion} to the most general cases. And we give the sharp embedding between Fourier $L^{p}$ spaces and modulation spaces.
	\end{abstract}
	
	\section{Introduction}\ 
	
	The modulation spaces $\mpq^{s}$ are one of the function spaces introduced by Feichtinger \cite{Feichtinger2003Modulation} in the 1980s using the short-time Fourier transform to measure the decay and the regularity of the function in a different way from the usual $L^{p}$ Sobolev spaces or Besov-Triebel spaces. The precise definitions will be given in Section \ref{sec pre}. Roughly speaking, the Besov-Triebel spaces mostly use the  dyadic decompositions of the frequency space while the modulation spaces using the uniform decompositions of the frequency space. Therefore, these spaces may have many different properties in some cases. Also, they have some properties similarly. So, the relationships between the modulation spaces and Besov-Triebel spaces are important.
	
	Many basic properties of the modulation spaces had been studies in \cite{Triebel1983Modulation,Groechenig2001Foundations,Kobayashi2006Modulation} such as the Banach spaces properties and the dual spaces of modulation spaces. And the modulation spaces also have some applications in pseudo-differential operators in \cite{Sugimoto1988Pseudo,Groechenig1999Modulation,Toft2001Subalgebras,Toft2002Modulation,Toft2002Continuity},\cite{Toft2004Continuity}. Modulation spaces also have many applications in the analysis of partial differential equations. For example, the Schrodinger and wave semigroups which are not bouned on neither $L^{p}$ or $\besov$ for $p\neq 2$, are bounded on $\mpq^{s}$(see \cite{Benyi2007Unimodular}). So, the modulation space is a good space for the initial data of the Cauchy problem for nonlinear disperse equations(see \cite{Baoxiang2006Isometric,Wang2007Frequency,Wang2009Global,Ruzhansky2012Modulation}) .  
	
	As for embedding properties of modulation spaces, Toft in \cite{Toft2001Embeddings} discussed  some of the sufficient conditions of the embedding between modulation spaces and Besov spaces, also, Toft in \cite{Toft2004Continuity} discussed the embedding between the $L^{p}$ spaces and the Fourier transform of $L^{p}$ space with the modulation spaces. Then, Okoudjou in \cite{Okoudjou2004Embedding} gave some sufficient conditions for some cases of embedding between Sobolev-Besov spaces and modulation spaces by the means of the short-time Fourier transform. Later, Sugmito and Tomita in \cite{Sugimoto2007dilation} gave the sharp embedding between modulation spaces and Besov spaces when $1\leq p,q \leq \infty$ by means of the scaling of modulation spaces, then Wang and Huang in \cite{Wang2007Frequency} extended the condition to $0<p,q \leq \infty$ by means of the uniform decomposition of frequency space. But they all discuss the case when the indexes of physics space and frequency space of modulation spaces and Besov spaces are the same. Later, Kobayashi and Sugimoto in \cite{Kobayashi2011inclusion} studied the sharp embedding between $L^{p}$-Sobolev space and modulation spaces using the method similar to \cite{Kobayashi2009Embedding}. As for Triebel spaces, Guo, Wu and Zhao in \cite{Guo2017Inclusion} gave the sharp condition of the embedding $F_{p,r} \hookrightarrow \mpq^{s}$ when $0<p \leq 1$.
	
	In this article, we study the sharp conditions of the embedding between Besov-Triebel-Sobolev spaces and modulation spaces. Let us take Besov spaces for example, in Theorem \ref{besov to mpq}, we get the sufficient and necessary conditions of the embedding $B_{p_{0},q_{0}}^{s} \hookrightarrow \mpq$, which is the most general case of $0<p,q,p_{0},q_{0} \leq \infty$, also in Theorem \ref{mpq to besov}, we get the sufficient and necessary conditions of the embedding of the opposite side. In the case of $0<p,q<1$, the dual of Besov spaces and modulation spaces are  different, so we cannot just get the opposite side by the method of duality. Instead, we use the method of the uniform decompositions of frequency space to get the result as desired. 
	
	We explain the organization of this paper. We give some basic notations and properties of the spaces we need in the next preliminary section. In section 3, we give some lemmas we need in the proof, especially for the case of $0<p,q\leq \infty$ we discuss in section 4. The embedding between Besov spaces and modulation spaces is considered in section 4. We mainly use the method of uniform decomposition in \cite{Wang2007Frequency}. Then we study the embedding between $L^{p}$ Sobolev spaces and modulation spaces in section 5, as we can see, there are four conditions totally. In section 6, we consider the embedding Triebel spaces and modulation spaces which is sometimes equivalent to the embedding of Sobolev spaces by Littlewood-Paley theory. In the last section, we give the sharp embedding relationship between the Fourier $L^{r}$ spaces and modulation spaces.

	\section{Perliminaries}\label{sec pre}
	\subsection{Basic natation}
	
	 The following notation will be used throughout this article. For $0< p,q \leq \infty$, we denote 
	\begin{align*}
		\sigma (p,q) &:=d\left(0 \wedge (\frac{1}{q} - \frac{1}{p}) \wedge (\frac{1}{q} +\frac{1}{p} -1)\right);\\
		\tau(p,q) & := d\left(0 \vee (\frac{1}{q} - \frac{1}{p}) \vee (\frac{1}{q} +\frac{1}{p} -1)\right),
	\end{align*}
	where $a\wedge b= \min\set{a,b},a\vee b =\max\set{a,b}$.

	We write $\sch(\real^{d})$ to denote the Schwartz space of all complex-valued rapidly decreasing infinity differentiable functions on $\real^{d}$, and $\sch'(\real^{d})$ to denote the dual space of $\sch(\real^{d})$, all called the space of all tempered distributions. For simplification, we omit $\real^{d}$ without causing ambiguity. The Fourier transform is defined by $\F f(\xi) = \hat{f}(\xi) = \int_{\real^{d}} f(x) e^{-ix\xi} d\xi$, and the inverse Fourier transform by $\FF f(x) = (2\pi)^{-d} \int_{\real^{d}}f(\xi) e^{ix\xi} d\xi$. For $1\leq p <\infty$, we define the $L^{p}$ norm:
	\begin{align*}
		\norm{f}_{p} = 	\left( \int_{\real^{d}} \abs{f(x)}^{p}\right)^{1/p}
	\end{align*}
	and $\norm{f}_{\infty} = ess\sup_{x\in \real_{d}} \abs{f(x)}$. We also define the $L^{p}$ Sobolev norm :
	\begin{align*}
		\norm{f}_{W^{s,p}}=\norm{(I-\triangle)^{s/2}f}_{p}.
	\end{align*}
	
	Recall that the Sobolev spaces is defined by $W^{s,p}=\set{f\in \sch': \norm{f}_{W^{s,p}} <\infty}$. And the Fourier $L^{p}$ is defined by $\F L^{p} =\set{ f\in \sch': \norm{\hat{f}}_{p} <\infty}$.

	We use the notation $I \lesssim J$ if there is an independently constant C such that $I \leq C J$. Also we denote $I \approx J$ if $I \lesssim J$ and $J\lesssim I$. For $1\leq p \leq \infty$, we denote $1/p+1/p'=1$, for $0<p<1$, denote $p'=\infty$.

	\subsection{Modulation spaces} \

	Let $0<p,q\leq \infty$, the short time Fourier tranform of $f$ respect to a window function $g\in \sch$ is defined as (see \cite{Feichtinger2003Modulation}):
	\begin{align*}
		V_{g}f(x,\xi) = \int_{\real^{d}} f(t) \overline{g(t-x)} e^{-it\xi} dt.
	\end{align*}
	And then for $1\leq p,q\leq \infty$, we denote 
	\begin{align*}
		\norm{f}_{\mpq^{s}} = \norm{V_{g}f(x,\xi)\inner{\xi}^{s}}_{L_{\xi}^{q} L_{x}^{p}},
	\end{align*}
	where $\inner{\xi} = 1+\abs{\xi}$.
	Modulation space $\mpq^{s}$ are defined as the space of all tempered distribution $f\in 
	\sch'$ for which $\norm{f}_{\mpq^{s}}$ is finite. 
	
	Also, we know another equivalent definition of modulation by uniform decomposition of frequency space (see \cite{Groechenig2001Foundations}).
	
	Let $\sigma$ be a smooth cut-off function adapted to the unit cube $[-1/2,1/2]^{d}$ and $\sigma=0$ outside the cube $[-3/4,3/4]^{d}$, we write $\sigma_{k} = \sigma(\cdot -k)$, and assume that 
	\begin{align*}
		\sum_{k\in\Z^{d}} \sigma_{k} (\xi) \equiv 1, \ \forall \xi \in \real^{d}.
	\end{align*}
	
	Denote $\sigma_{k} (\xi) =\sigma (\xi-k)$, and $\Box_{k} = \FF \sigma_{k} \F $, then we have the following equivalent norm of modulation space:
	
	\begin{align*}
		\norm{f}_{\mpq^{s}} = \norm{\inner{k}^{s}\norm{\Box_{k}f}_{L_{x}^p}}_{\ell_{k}^{q} (\Z^{d})}.
	\end{align*}
	
	We simply write $\mpq$ instead of $\mpq^{0}$. One can prove the $\mpq^{s}$ norm is independent of the choice of cut-off function $\sigma$. Also $\mpq^{s}$ is a quasi Banach space and when $1\leq p,q \leq \infty$, $\mpq^{s}$ is a Banach space. When $p,q <\infty$, then $\sch$ is dense in $\mpq^{s}$. Also, $\mpq^{s}$ has some basic properties, we list them in the following lemma(see \cite{Wang2007Frequency,Wang2011Harmonic}).
	\begin{lemma}\label{lemma mpq}
		For $s,s_{0},s_{1},\in \real,0<p,p_{0},p_{1},q,q_{0},q_{1} \leq \infty$,
		\begin{description}
			\item[(1)] if $s_{0} \leq s_{1},p_{1} \leq p_{0},q_{1}\leq q_{0}$, we have $M_{p_{1},q_{1}}^{s_{1}} \hookrightarrow M_{p_{0},q_{0}}^{s_{0}}$;
			\item[(2)] when $p,q< \infty$, the dual space of $\mpq^{s}$ is $M_{p',q'}^{-s}$;
			\item[(3)] the interpolation spaces theorem is true for $\mpq^{s}$, i.e. for $0<\theta<1$
			when $$
			s=(1-\theta) s_{0}+\theta s_{1}, \quad \frac{1}{p}=\frac{1-\theta}{p_{0}}+\frac{\theta}{p_{1}}, \quad \frac{1}{q}=\frac{1-\theta}{q_{0}}+\frac{\theta}{q_{1}},$$ we have $\left(M_{p_{0}, q_{0}}^{s_{0}}, M_{p_{1}, q_{1}}^{s_{1}}\right)_{\theta}=M_{p, q}^{s}$;
			
			\item[(4)] when $q_{1}<q,s+d/q>s_{1} +d/q_{1}$, then we have $\mpq^{s} \hookrightarrow M_{p,q_{1}}^{s_{1}}$.
		\end{description}
	\end{lemma}

	\subsection{Besov-Triebel spaces}\ 
	
	Let $0<p,q \leq \infty,s\in \real,$ choose $\psi: \real^{d}:\rightarrow [0,1] $ be a smooth radial bump function adapted to the ball $B(0,2)$: $\psi(\xi) =1$ as $\abs{ \xi} \leq1$ and $\psi(\xi) =0$ as $\abs{\xi} \geq2$. We denote $\varphi(\xi) = \psi(\xi)- \psi(2\xi)$, and $\varphi_{j}(\xi) = \varphi(2^{-j} \xi)$ for $1\leq j,j\in \Z $, $\varphi_{0}(\xi ) = 1- \sum_{j\geq1} \varphi_{j}(\xi)$. We say that $\triangle_{j} = \FF  \varphi_{j} \F$ are the dyadic decomposition operators. The Besov spaces $\besov$ and the Triebel spaces $\triebel$ are definied in the following way :
	\begin{align*}
		\besov = \set{ f\in \sch'(\real^{d}) : \norm{f}_{\besov} = \norm{ 2^{js} \norm{\triangle_{j}f}_{L_{x}^p} }_{\ell_{j\geq0}^{q}} < \infty},\\
		\triebel = \set{ f\in \sch'(\real^{d}) : \norm{f}_{\triebel} = \norm{\norm{2^{js} \triangle_{j} f} _{\ell_{j\geq0} ^{q}}}_{L_{x}^{p}} < \infty}.
	\end{align*}
	
	One can prove that the Besov-Triebel norms defined by different dyadic decompositions are all equivalent(see \cite{Triebel1983Theory}), so without loss of generality, we can assume that when $1\leq j$, $\varphi_{j}(\xi)=1$ on $D_{j} := \set{\xi\in \real^{d}: \frac{3}{4} 2^{j} \leq \abs{\xi} \leq \frac{5}{4} 2^{j}}$ for convenience. Also, Besov-Triebel spaces have some basic properties known already(see \cite{Triebel1983Theory,Kobayashi2006Modulation,Wang2011Harmonic}).
	\begin{lemma}
		\label{lemma besov triebel}
		Let $s,s_{1},s_{2} \in \real,0<p,p_{1},p_{2},q,q_{1},q_{2} \leq \infty$,
		\begin{description}
			\item[(1)] if $q_{1} \leq q_{2}$, we have $B_{p,q_{1}}^{s} \hookrightarrow B_{p,q_{2}}^{s}$, $F_{p,q_{1}}^{s} \hookrightarrow F_{p,q_{2}}^{s}$;
			\item[(2)] $\forall \eps >0$, we have $B_{p,q_{1}}^{s+\eps} \hookrightarrow B_{p,q_{2}}^{s}$, $F_{p,q_{1}}^{s+\eps} \hookrightarrow F_{p,q_{2}}^{s}$;
			\item[(3)] $B_{p,p\wedge q}^{s} \hookrightarrow \triebel \hookrightarrow B_{p,p\vee q}^{s}$;
			\item[(4)] if $p_{1} \leq p_{2},s_{1}-d/p_{1} =s_{2}-d/p_{2}$, we have $B_{p_{1},q}^{s_{1}} \hookrightarrow B_{p_{2},q}^{s_{2}}$;
			\item[(5)] if $p_{1}<p_{2},s_{1}-d/p_{1} =s_{2}-d/p_{2}$, we have $F_{p_{1},q_{1}}^{s_{1}} \hookrightarrow F_{p_{2},q_{2}}^{s_{2}}$;
			\item[(6)] when $1\leq p,q <\infty$, the dual space of $\besov$ is $B_{p',q'}^{-s}$, the dual space of $\triebel$ is $F_{p',q'}^{-s}$;
			\item[(7)]  the interpolation spaces theorem is true for $\besov$ and $\triebel$, i.e. for $0<\theta<1$ when $$s=(1-\theta) s_{0}+\theta s_{1}, \quad \frac{1}{p}=\frac{1-\theta}{p_{0}}+\frac{\theta}{p_{1}}, \quad \frac{1}{q}=\frac{1-\theta}{q_{0}}+\frac{\theta}{q_{1}},$$ we have $\left(B_{p_{0}, q_{0}}^{s_{0}}, B_{p_{1}, q_{1}}^{s_{1}}\right)_{\theta}=B_{p, q}^{s}$, $\left(F_{p_{0}, q_{0}}^{s_{0}}, F_{p_{1}, q_{1}}^{s_{1}}\right)_{\theta}=F_{p, q}^{s}$ ;
			\item[(8)] when $1<p<\infty$, $F_{p,2}^{s} = W^{s,p}$, $F_{1,2}^{s} \hookrightarrow W^{s,1},W^{s,\infty} \hookrightarrow F_{\infty,2}^{s}$.
		\end{description}
	\end{lemma}

	\section{Some lemmas}\ 
	
	In this section, we give some useful lemmas and propositions proved already, which will be used in our proof.

	The following Bernstein's inequality is very useful in time-frequency analysis (see \cite{Wang2011Harmonic}) :
	
	\begin{lemma}\label{lemm bernstein}
		[Bernstein's inequality]
		Let $0< p \leq q \leq \infty, b>0,\xi_{0} \in \real^{d}$. Denote $L_{B(\xi,b)}^{p} = \set{ f\in L^{p}: \mbox{ supp } \hat{f} \subseteq B(\xi,R)}$. Then there exists $C(d,p,q)>0$, such that \begin{align*}
			\norm{f}_{q} \leq C(d,p,q) R^{d(1/p-1/q)} \norm{f}_{p}
		\end{align*} 
		holds for all $f\in L_{B(\xi,b)}^{p}$ and $C(d,p,q)$ is independent of $b>0$ and $\xi_{0} \in \real^{d}$.
	\end{lemma}

	Also, by using the Bernstein's inequality, we can get the following Young type inequality for $0<p<1$:
	
	\begin{lemma}[\cite{Kobayashi2006Modulation}]
		\label{lemma young}
		Let $0<p<1,R_{1},R_{2}>0,\xi_{1},\xi_{2} \in \real^{d}$, then there exists $C(d,p)>0$, such that 
		\begin{align*}
			\norm{\abs{f} * \abs{g}}_{p} \leq C(d,p) (R_{1}+R_{2})^{d(1/p-1)} \norm{f} \norm{g}
		\end{align*}
		holds for all $f\in L_{B(\xi_{1},R_{1})}^{p}, g\in L_{B(\xi_{2},R_{2})}^{p}$.
	\end{lemma}
	
	By the definition of $\Box_{k}$ and $\triangle_{j}$ operators, we know that they are all convolution type operators, so by the Lemma above and the usual Young's inequality, we have:
	
	\begin{coro}
		\label{cor convolution}
		For $0<p\leq \infty,1\leq j,\ell,k\in \Z^{d}$, there exists $C>0$ independent of $k,j$, such that 
		\begin{align*}
			\norm{\triangle_{j} \Box_{k} f} &\leq C \norm{\Box_{k}f}_{p},\\
			\norm{\triangle_{j}\triangle_{\ell} f} _{p} &\leq C \norm{\triangle_{j} f}_{p}.
		\end{align*}
	\end{coro}

	Finally we recall a result known already, which is the special case of Theorem \ref{besov to mpq} and Theorem \ref{mpq to besov}.
	
	\begin{prop}[\cite{Wang2007Frequency}]
		\label{prop besov to mpq}
		For $0< p,q \leq \infty,s\in \real$, we have 
		\begin{description}
			\item[(1)] $\besov \hookrightarrow \mpq$ if and only if $s\geq \tau(p,q)$;
			\item[(2)] $\mpq \hookrightarrow \besov$ if and only if $s\leq \sigma(p,q)$.
		\end{description}
	\end{prop}
	
	\begin{prop}(\cite{Kobayashi2011inclusion})
		\label{prop wsp to mpq}
		Let $1\leq p,q \leq \infty,s\in \real$, then $W^{s,p} \hookrightarrow \mpq$ if and only if one of the following cases is satisfied:
		\begin{description}
			\item[(1)] $1<p \leq q, s\geq \tau(p,q)$;
			\item[(2)] $p>q, s>\tau (p,q)$;
			\item[(3)] $p=1,q=\infty, s\geq \tau(1,\infty)$;
			\item[(4)] $p=1,q\neq \infty, s>\tau(1,q)$.
		\end{description}
	\end{prop}

	\section{Embedding between Besov spaces and modulation spaces} \ 
	
	In this section, we consider the embedding with the form $\besovo \hookrightarrow \mpq$. If we got the sufficient and necessary condition of this embedding, then the embedding $B_{p_{1},q_{1}}^{s_{1}} \hookrightarrow \mpq^{s}$ can be done in the same way.
	
	Our main result is:
	\begin{thm}
		\label{besov to mpq}
		For $0< p,q,p_{0},q_{0} \leq \infty,s\in \real$, the embedding $\besovo \hookrightarrow \mpq$ is true if and only if one of the following two conditions is satisfied:
		\begin{description}
			\item[(1)] $p_{0} \leq p, q_{0} \leq q, s\geq \tau(p_{0},q)$;
			\item[(2)] $p_{0} \leq p, q_{0} >q, s> \tau(p_{0},q)$.
		\end{description}
	\end{thm}
	
	As for the other side of the embedding between Besov spaces and modulation spaces, we also have:
	\begin{thm}
		\label{mpq to besov}
		For $0< p,q,p_{1},q_{1} \leq \infty,s\in \real$, the embedding $\mpq \hookrightarrow \besovon$ is true if and only if one of the following two conditions is satisfied:
		\begin{description}
			\item[(1)] $p_{1} \geq p, q_{1} \geq q, s\leq \sigma(p_{1},q)$;
			\item[(2)] $p_{1} \geq p, q_{1} <q, s < \sigma(p_{1},q)$.
		\end{description}	
	\end{thm}
\begin{rem}
If we regard the Besov space as a special $\alpha$-modulation space when $\alpha=1$ as stated in \cite{Guo2018Full}, the result above is coincident with the Theorem 1.2 in \cite{Guo2018Full}, where the authors got the full characterization of embedding between $\alpha$-modulation spaces. 
\end{rem}

	For a special case when $p_{0}=q_{0}=2$, which means that $\besovo = H^{s}$, the classical Sobolev spaces, then we have 
	\begin{coro}
		Let $0< p,q \leq \infty, s\in \real$, we have 
		
		\begin{description}
			\item[(1)] 	 the embedding $H^{s} \hookrightarrow \mpq$ is true if and only if one of the following two cases is true:
			\begin{enumerate}
				\item $2\leq p, 2\leq q, s\geq 0$;
				\item $ 2\leq p, 2>q, s>d(1/q-1/2)$.
			\end{enumerate}
			\item[(2)] The embedding $\mpq \hookrightarrow H^{s}$ is true if and only if the following two cases is true:
			\begin{enumerate}
				\item $p\leq 2, q\leq 2, s\leq 0$;
				\item $p\leq 2, q>2, s<d(1/q-1/2)$.
			\end{enumerate}
		\end{description}
		
	\end{coro}
	
	\subsection{Proof of Theorem \ref{besov to mpq}}\ 
	
	In this subsection, we prove the Theorem \ref{besov to mpq}. For sufficiency, by using the lemmas above, we can easily prove it. As for necessity, it may be much more difficult. Here, we mostly use the way in [\cite{Wang2007Frequency}] to prove it.
	
	\subsubsection{Sufficiency}
	
	Use Lemma \ref{lemma besov triebel}, we know that for any case in \ref{besov to mpq}, we have $\besovo \hookrightarrow B_{p_{0},q}^{\tau(p_{0},q)}$; then by Proposition \ref{prop besov to mpq}, we know that $B_{p_{0},q}^{\tau(p_{0},q)} \hookrightarrow M_{p_{0},q}$; also by Lemma \ref{lemma mpq}, we have $M_{p_{0},q} \hookrightarrow \mpq$.

	\subsubsection{Necessity}\label{nece}
	If we have $\besovo \hookrightarrow \mpq$, then we have
	\begin{align}\label{eq besov to mpq}
		\norm{f}_{\mpq} \lesssim \norm{f}_{\besovo}, \  \forall f\in \besovo.
	\end{align}
	Then we can take some different kinds of $f$ into \eqref{eq besov to mpq} to get some restrictions of indexes.
	
	\begin{description}
		\item[(1)]  Choose $\eta \in \sch$, such that supp $\eta \subset [-1/8,1/8]^{d}$, denote $f= \FF \eta$, for $\forall \lambda \leq1$, take $f_{\lambda}(x) = f(\lambda x)$ into \eqref{eq besov to mpq}. By the choose of $\eta$, we know that supp $\F f_{\lambda} \subset \lambda [-1/8,1/8]^{d}$, which means that 
		\begin{align*}
			\Box_{k}f_{\lambda} = \begin{cases*}
				f_{\lambda}, & $k=0$;\\
				0, &$k\neq 0$;
			\end{cases*}
			\triangle_{j}f_{\lambda} = \begin{cases*}
				f_{\lambda}, & $k=0$;\\
				0, &$k\neq 0$.
			\end{cases*}
		\end{align*} 
		So, take $f_{\lambda}$ into \eqref{eq besov to mpq}, we have $\norm{f_{\lambda}}_{p} \lesssim \norm{f_{\lambda}}_{p_{0}}$, which means that $\lambda^{-d/p} \lesssim \lambda^{-d/p_{0}}$. Take $\lambda \rightarrow 0$, we have $p_{0} \leq p$.
		
		\item[(2)] Denote $A_{j}:= \set{k\in \Z^{d}: k+[-3/4,3/4]^{d} \subseteq D_{j}}$, for $\ell \geq 10$, choose $k_{\ell} \in A_{\ell}$, take $f_{\ell} = \FF \eta(\cdot -k_{\ell})$, $\eta$ is just the function in (1). Then we know that supp $\F f_{\ell} \subseteq k_{\ell} + [-1/8,1/8]^{d} \subseteq D_{\ell}$, so we have 
		\begin{align*}
			\Box_{k} f_{\ell} = \begin{cases*}
				f_{\ell}, & $k=k_{\ell}$;\\
				0, & $k\neq k_{\ell};$
			\end{cases*}
			\triangle_{j}f_{\ell} = \begin{cases*}
				f_{\ell}, & $j=\ell$;\\
				\triangle_{j}f_{\ell}, & $ 0<\norm{j-\ell}_{\infty} \leq 3$;\\
				0, & $\norm{j-\ell}_{\infty} >3.$
			\end{cases*}
		\end{align*}
		Then by Corollary \ref{cor convolution} ,we have 
		\begin{align*}
			\norm{f_{\ell}}_{\mpq} &= \norm{f_{\ell}}_{p} \approx 1\\
			\norm{f_{\ell}}_{\besovo} & =\norm{2^{js} \norm{\triangle_{j} f_{\ell}}_{p_{0}}}_{\ell_{\norm{j-\ell}_{\infty} \leq 3}^{q_{0}}}\\
			& \lesssim 2^{\ell s} \norm{f_{\ell}}_{p_{0}} \approx 2^{\ell s}.
		\end{align*}
		Then take the estimate above into \eqref{eq besov to mpq}, we have $1\lesssim 2^{\ell s}, \forall \ell \geq 10$, so we get $s\geq 0$.
		
		\item[(3)] Take $f_{\ell} = \FF \varphi_{\ell}$, then by Corollary \ref{cor convolution}, we have 
		\begin{align*}
			\norm{f_{\ell}}_{\mpq} &= \norm{\norm{\Box_{k}f_{\ell}}_{p}}_{\ell_{k}^{q}} \geq \norm{\norm{\Box_{k} f_{\ell}}_{p}}_{\ell_{k\in A_{\ell}}^{q}} 
			= \norm{ \norm{\FF \sigma_{k}}_{p}}_{\ell_{k\in A_{\ell}}^{q}}  \approx 2^{\ell d/q};\\
			\norm{f_{\ell}}_{\besovo} &= \norm{ 2^{js} \norm{\triangle_{j}f_{\ell}}_{p_{0}}} _{\ell_{j}^{q_{0}}} = \norm{ 2^{js}\norm{\triangle_{j}f_{\ell}}_{p_{0}}}_{\ell_{\norm{j-\ell}_{\infty}\leq 3}^{q_{0}}} \leq 2^{\ell s} \norm{f_{\ell}}_{p_{0}} \approx 2^{\ell(s+d(1-1/p_{0}))}.
		\end{align*}
		Take the estimates above into \eqref{eq besov to mpq}, let $\ell \rightarrow \infty$, we have $s\geq d(1/p_{0}+1/q-1)$.

		\item[(4)] When $p_{0}\geq 2$, choose $\eta \in \sch$, such that supp $\hat{\eta} \subseteq[-1/8,1/8]^{d}$, take $f_{\ell}(x) = \sum_{k\in A_{\ell}} e^{ikx} \eta(x-k)$, then $\hat{f}_{\ell}(\xi) = \sum_{k\in A_{\ell}} e^{-ik (\xi-k)} \hat{\eta}(\xi-k)$, so we know that supp $\hat{f}_{\ell} \subseteq \bigcup_{k\in A_{\ell}} k+[-1/8,1/8]^{d} \subseteq D_{\ell}$, so by the orthogonality of $\triangle_{j},\Box_{k}$, we have:
		\begin{align*}
			\Box_{k}f_{\ell} = \begin{cases*}
				e^{ikx} \eta(x-k), & $k\in A_{\ell}$;\\
				0, & else;
			\end{cases*}\ \ 
			\triangle_{j}f_{\ell} = \begin{cases*}
				f_{\ell}, & $j=\ell$;\\
				\triangle_{j}f_{\ell}, & $ 0<\norm{j-\ell}_{\infty} \leq 3$;\\
				0, & $\norm{j-\ell}_{\infty} >3.$
			\end{cases*}
		\end{align*}
		Then we have \begin{align*}
			&\norm{f_{\ell}}_{\mpq} \geq \norm{\norm{\Box_{k}f_{\ell}}_{p}}_{\ell_{k\in A_{\ell}}^{q}} \approx 2^{\ell d/q};\\
			&\norm{f_{\ell}}_{\besovo} = \norm{ 2^{js} \norm{\triangle_{j}f_{\ell}}_{p_{0}}} _{\ell_{j}^{q_{0}}} = \norm{ 2^{js}\norm{\triangle_{j}f_{\ell}}_{p_{0}}}_{\ell_{\norm{j-\ell}_{\infty}\leq 3}^{q_{0}}} \leq 2^{\ell s} \norm{f_{\ell}}_{p_{0}} \leq 2^{\ell s} \norm{f_{\ell}}_{2}^{2/p_{0}} \norm{f}_{\infty}^{1-2/p_{0}}.
		\end{align*}
		By the orthogonality of $L^{2}$ space, we have \begin{align*}
			\norm{f_{\ell}}_{2} = \left(\sum_{k\in A_{\ell}} \norm{\eta(x-k)}_{2}^{2}  \right)^{1/2} \approx 2^{\ell d/2}.
		\end{align*}
		Also by the fast decay of $\eta$, we have \begin{align*}
			\abs{f_{\ell}(x)} \leq \sum_{k\in A_{\ell}} \abs{\eta(x-k)} \leq \sum_{k\in\Z^{d}} (1+\abs{x-k})^{-N} \lesssim 1.
		\end{align*}
		Combine the estimates above, we have $ \norm{f}_{\besovo} \lesssim 2^{\ell(s+d/p_{0})}$. Then take it into \eqref{eq besov to mpq}, and put $\ell \rightarrow \infty$, we have $s\geq d(1/q-1/p_{0})$.

		\item[(5)] Combine case(2)-(4), we have $s\geq \tau(p_{0},q)$. In contrast with the necessary condition we need in Theorem \ref{besov to mpq}, we only need to prove that if $B_{p_{0},q_{0}}^{\tau(p_{0},q)} \hookrightarrow \mpq$, then we have $q_{0}\leq q$.
		\begin{description}
			\item[(5.1)] When $\tau(p_{0},q) =0$, take $f= \sum_{\ell} a_{\ell} f_{\ell}$, where $f_{\ell}$ in case (2), then by the same calculation, we have \begin{align*}
				\norm{f}_{\mpq} \approx \norm{a_{\ell}}_{\ell^{q}};\ \ \norm{f}_{B_{p_{0},q_{0}}^{0}} \lesssim \norm{a_{\ell}}_{\ell^{q_{0}}}. 
			\end{align*}
			Take the estimates into \eqref{eq besov to mpq}, we have $q_{0}\leq q$.
			\item[(5.2)] When $\tau(p_{0},q) =d(1/p_{0}+1/q-1)$, take $f= \sum_{\ell} a_{\ell} f_{\ell}$, where $f_{\ell}$ in case (3), then by the same calculation, we have
			\begin{align*}
				\norm{f}_{\mpq} &= \norm{\norm{\sum_{\ell} a_{\ell} \Box_{k}f_{\ell}}_{p}}_{\ell_{k}^{q}} \geq \norm{\norm{\norm{\sum_{\ell} a_{\ell} \Box_{k}f_{\ell}}_{p}}_{\ell_{k\in A_{j}}^{q}}}_{\ell_{j}^{q}}= \norm{a_{j}2^{jd/q}}_{\ell_{j}^{q}};\\
				\norm{f}_{B_{p_{0},q_{0}}^{\tau(p_{0},q)}} 
				&= \norm{ 2^{j\tau} \norm{\sum_{\ell} a_{\ell} \triangle_{j} f_{\ell}}_{p_{0}}}_{\ell_{j}^{q_{0}}} \leq \sum_{\norm{\ell}_{\infty}\leq 3}  \norm{a_{j+\ell} 2^{j\tau} \norm{\triangle_{j} f_{\ell+j}}_{p_{0}}}_{\ell_{j}^{q_{0}}}\\
				& \leq \sum_{\norm{\ell}_{\infty}\leq 3}  \norm{a_{\ell+j} 2^{jd(1-1/p_{0})} 2^{j\tau}}_{\ell_{j}^{q_{0}}} \lesssim \norm{a_{j}2^{jd/q}}_{\ell_{j}^{q_{0}}}.
			\end{align*}
			Take the estimates into \eqref{eq besov to mpq}, we have $q_{0}\leq q$.
			\item[(5.3)] When $\tau(p_{0},q)=d(1/q-1/p_{0}),$ take $f= \sum_{\ell} a_{\ell} f_{\ell}$, where $f_{\ell}$ in case (4), then by the same calculation, we have
			\begin{align*}
				\norm{f}_{\mpq} &= \norm{\norm{\sum_{\ell} a_{\ell} \Box_{k}f_{\ell}}_{p}}_{\ell_{k}^{q}} \geq \norm{\norm{\norm{\sum_{\ell} a_{\ell} \Box_{k}f_{\ell}}_{p}}_{\ell_{k\in A_{j}}^{q}}}_{\ell_{j}^{q}}\\
				&= \norm{ \norm{a_{j} \norm{\FF \sigma_{k}}_{p}}_{\ell_{k\in A_{j}}^{q}}}_{\ell_{j}^{q}} \approx \norm{a_{j}2^{jd/q}}_{\ell_{j}^{q}};\\
				\norm{f}_{B_{p_{0},q_{0}}^{\tau}} &= \norm{ 2^{j\tau} \norm{\sum_{\ell} a_{\ell} \triangle_{j} f_{\ell}}_{p_{0}}}_{\ell_{j}^{q_{0}}} \leq \sum_{\norm{\ell}_{\infty}\leq 3}  \norm{a_{j+\ell} 2^{j\tau} \norm{\triangle_{j} f_{\ell+j}}_{p_{0}}}_{\ell_{j}^{q_{0}}}\\
				&\leq \sum_{\norm{\ell}_{\infty}\leq 3}  \norm{a_{\ell+j} \norm{f_{\ell+j}}_{p_{0}} 2^{j\tau}}_{\ell_{j}^{q_{0}}} \lesssim \norm{a_{j}  2^{j\tau}2^{jd/p_{0}}}_{\ell_{j}^{q_{0}}} = \norm{a_{j}2^{jd/q}}_{\ell_{j}^{q_{0}}}.
			\end{align*}
			Take the estimates into \eqref{eq besov to mpq}, we have $q_{0}\leq q$.
		\end{description}
		
	\end{description}

	\subsection{Proof of Theorem \ref{mpq to besov}}\ 
	
	In this section, we prove the Theorem \ref{mpq to besov}. Because we consider the case when $0<p,q\leq \infty$, we can not only use the duality to prove it. Thanks to the wonderful method in [\cite{Wang2007Frequency}], also combining the construction in the proof of Theorem \ref{besov to mpq}, we can prove this theorem as well. 
	
	\subsubsection{Sufficiency}
	By Lemma \ref{lemma mpq}, we have $\mpq \hookrightarrow M_{p_{1},q}$; by Proposition \ref{prop besov to mpq}, we have $M_{p_{1},q} \hookrightarrow B_{p_{1},q}^{\sigma(p_{1},q)}$; finally by Lemma \ref{lemma besov triebel}, we have $B_{p_{1},q}^{\sigma(p_{1},q)} \hookrightarrow \besovon$ for any case in Theorem \ref{mpq to besov}.

	\subsubsection{Necessity}
	If we have $\mpq \hookrightarrow \besovon$, then we have 
	\begin{align}\label{eq mpq to besov}
		\norm{f}_{\besovon} \lesssim \norm{f}_{\mpq}.
	\end{align}
	
	Just like the proof of Theorem \ref{besov to mpq}, we can take some different kinds of $f$ into \eqref{eq mpq to besov} to get some restrictions of indexes, and in some cases, we can take the same $f$ as in \ref{nece}.
	
	\begin{description}
		\item[(1)] Take $f_{\lambda}$ in the case (1) in \ref{nece} into \eqref{eq mpq to besov}, we can get $p\leq p_{1}$.
		
		\item[(2)] Take $f_{\ell}$ in the case (2) in \ref{nece}, we have $\norm{f_{\ell}}_{\mpq} \approx 1$, and $\norm{f_{\ell}}_{\besovon} \geq 2^{\ell s} \norm{\triangle_{\ell} f_{\ell}}_{p_{1}} = 2^{\ell s} \norm{f_{\ell}}_{p_{1}} \approx 2^{\ell s}$. Then we have $s\leq 0$. 
		
		\item[(3)] Denote $B_{j} = \set{k\in \Z^{d}: k+[-3/4,3/4]^{d} \cap D_{j} \neq \emptyset}$, choose $g\in \sch,$ such that supp $g \subseteq D_{0}$,  then denote $g_{\ell}(\cdot) = g(2^{-\ell} \cdot)$. Take $f_{\ell} = \FF g_{\ell}$, so supp $\F f_{\ell} \subseteq D_{\ell}$. So by Corollary \ref{cor convolution}, we have 
		\begin{align*}
			&\norm{f_{\ell}}_{\mpq} = \norm{\norm{\Box_{k}f_{\ell}}_{p}}_{\ell_{k\in B_{\ell}}^{q}} \leq \norm{\norm{\FF \sigma_{k}}_{p}}_{\ell_{k\in B_{\ell}}^{q}} \lesssim 2^{\ell d/q};\\
			&\norm{f_{\ell}}_{\besovon} \geq 2^{\ell s} \norm{\triangle_{\ell} f_{\ell}}_{p_{1}}=2^{\ell s} \norm{f_{\ell}}_{p_{1}} \approx 2^{\ell (s+d(1-1/p_{1}))}.
		\end{align*}
		Take the estimate into \eqref{eq mpq to besov}, let $\ell \rightarrow \infty$, we can get $s\leq d(1/p_{1}+1/q-1)$.
		
		\item[(4)] Take $f_{\ell} (x) = \sum_{k\in A_{\ell}} e^{ikx} \eta(\frac{x-k}{a})$, we can choose $a<<1$ such that $\norm{f_{\ell}}_{p_{1}} \geq C 2^{\ell d/p_{1}}$(see \cite{Wang2007Frequency}), then by the same calculation of case (4) in \ref{nece}, we have \begin{align*}
			\norm{f_{\ell}}_{\besovon} \geq C 2^{\ell (s+d/p_{1})};\ \ 
			\norm{f_{\ell}}_{\mpq} \lesssim 2^{\ell d/q}.
		\end{align*}
		Take the estimates above into \eqref{eq mpq to besov} and put $\ell \rightarrow \infty$, we have $s\leq d(1/q-1/p_{1})$.
		
		\item[(5)] Combine case(2)-(4), we have $s \leq \sigma(p_{1},q)$. In contrast with the necessary condition we need in Theorem \ref{mpq to besov}, we only need to prove that if $\mpq \hookrightarrow B_{p_{1},q_{1}}^{\sigma(p_{1},q)}$, then we have $q_{1} \geq q$.
		
		By the same way in the proof of Theorem \ref{besov to mpq}, take $f=\sum a_{\ell}f_{\ell}$, $f_{\ell}$ in case (2)-(4), we can get $\norm{a_{\ell}}_{q_{1}} \lesssim \norm{a_{\ell}}_{q}$, which means that $q_{1} \geq q$. 
		
	\end{description}

	\section{Embedding between Sobolev spaces and modulation spaces}\ 
	
	In this section, we can get the sufficient and necessary conditions of the embedding between $L^{r}$-Sobolev spaces and modulation spaces. We consider the most general case of this kind of embedding with the form like $W^{s,r} \hookrightarrow \mpq$, which is the general case of Proposition \ref{prop wsp to mpq} where $r=p$. Also, we consider the case when $0<q \leq \infty$. Unlike the proof in that proposition, we mostly use the uniform decomposition of frequency space and the Theorem \ref{besov to mpq} we get in last section. 
	
	Our main result is:
	\begin{thm}
		\label{thm wsr to mpq}
		Let $1\leq p,r\leq \infty,0<q\leq 
		\infty,s\in \real$, then $W^{s,r} \hookrightarrow \mpq$ if and only if $r\leq p$ and one of the following conditions is satisfied:
		\begin{description}
			\item[(1)] $ r>q , s> \tau(r,q)$;
			\item[(2)] $ 1<r \leq q, s\geq \tau(r,q)$;
			\item[(3)] $r=1,q=\infty, s\geq \tau(r,q)$;
			\item[(4)] $r=1,0<q<\infty,s> \tau(r,q)$;
		\end{description}	
	\end{thm}
	
	\begin{figure}
		\centering
		\begin{tikzpicture}[scale=5]
			\draw (0,0)--(1,1);
			\draw [->] (0,0) -- (1.2,0);
			\draw[->] (0,0) -- (0,1.52);
			\draw (1,1)--(1,0); 
			\filldraw[fill=gray!20] (0,0)--(1,1)--(1,1.5)--(0,1.5);
			\filldraw[fill =gray!40] (0,0)--(1,1)--(1,0);
						\draw[dashed] (0,1)--(1,1);
			\draw (1,1)--(1,1.5);
			\draw[- >] (0.5,0.5) -- (0.6,0.4) ;
			\draw[->] (1,0.5) -- (1.1,0.5);
			\draw[->] (1,0) -- (1.1,0.1);
			\draw (1/3,2/3) node {(1)} ;
			\draw (0.7,0.3) node {(2)};
			\draw (1.15,0.15) node {(3)};
			\draw (1.15,0.5) node {(4)}; 
			\draw (1.2,-2pt) node {$\frac{1}{r}$};
			\draw (-2pt,1.52) node {$\frac{1}{q}$};
			\draw (-1pt,-1pt) node {0};
			\draw (1,-2pt) node {1};
			\draw (-1pt,1) node {1};
		\end{tikzpicture}
		\caption{the four cases of $(1/r,1/q)$ in Theorem \ref{thm wsr to mpq}}
	\end{figure}
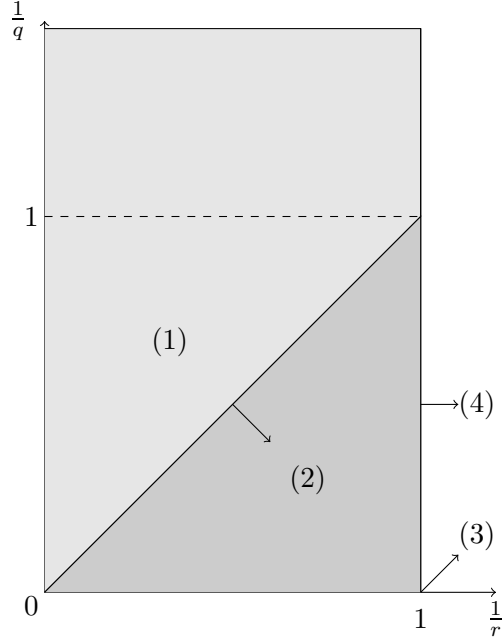

	\begin{rem}
		As the can see, the sufficient and necessary conditions we got in the theorem above is the same as the result in Proposition \ref{prop besov to mpq}. The index $p$ here only make  sense in $r\leq p$. The reason of this may be that the modulation space have good embedding property of the index $p$ like in Lemma \ref{lemma mpq}.
	\end{rem}

	As for the other side of this kind of embedding, we also have:
	\begin{thm}
		\label{thm mpq to wsr}
		Let $1\leq p,r \leq \infty, 0<q \leq \infty, s\in \real$, then $\mpq \hookrightarrow W^{s,r}$ if and only if $p\leq r$ and one of the following conditions is satisfied:
		\begin{description}
			\item[(1)] $r<q, s<\sigma(r,q)$;
			\item[(2)] $q\leq r<\infty, s\leq \sigma(r,q)$;
			\item[(3)] $r=\infty,0< q=1, s \leq \sigma(r,q)$;
			\item[(4)] $r = \infty, 1<q\leq \infty, s< \sigma(r,q)$;			
		\end{description}
	\end{thm}
	
	\subsection{Proof of Theorem \ref{thm wsr to mpq}}\ 
	
	In this subsection, we prove Theorem \ref{thm wsr to mpq}. By the Proposition \ref{prop wsp to mpq}, we can easily prove the sufficiency. As for necessity, we mainly construct some funcions to the embedding estimates to get some restrictions of the indexes.

	\subsubsection{Sufficiency}\ 
	
	When $1\leq q \leq \infty$, by Proposition \ref{prop wsp to mpq}, in any case of (1)-(4), we have $W^{s,r} \hookrightarrow M_{r,q}$; so, if $r\leq p$, by Lemma \ref{lemma mpq}, we have $M_{r,q} \hookrightarrow\mpq$.
	
	When $0<q<1$, which could only happen in case (2) or (4), then we have $s>\tau(r,q)=\tau(r,1)-d+d/q$. Take $\eps>0$ small enough, such that $s>\tau(r,1)-d+d/q+\eps$, then by Proposition \ref{prop wsp to mpq}, we have $W^{s,r} \hookrightarrow M_{r,1}^{d/q-d+\eps}$. Then by (4) in Lemma \ref{lemma mpq}, we have $M_{r,1}^{d/q-d+\eps} \hookrightarrow \mpq$, as desired.

	\subsubsection{Necessity} \label{subsubsec Wrq to mpq}
	If we have $W^{s,r} \hookrightarrow \mpq$, then we have 
	\begin{align} \label{eq wsr to mpq}
		\norm{f}_{\mpq} \lesssim \norm{f}_{W^{s,r}}.
	\end{align}
	\begin{description}
		\item[(i)] Choose $\eta \in \sch$, such that supp $\eta \subset [-1/8,1/8]^{d}$, denote $f= \FF \eta$, for $\forall \lambda \leq 1$, take $f_{\lambda}(x) = f(\lambda x)$ into \eqref{eq wsr to mpq}, we have $\norm{f_{\lambda}}_{r} \approx \norm{f_{\lambda}}_{W^{s,r}} \lesssim \norm{f_{\lambda}}_{\mpq} \approx \norm{f_{\lambda}}_{p}$, so we have $r\leq p$. Therefore, when $r=\infty$, we know that $p=\infty$, this is the case in Proposition \ref{prop wsp to mpq}.
		
		\item[(ii)] When $1\leq r <\infty$, by Lemma \ref{lemma besov triebel}, we have $B_{r,r\wedge 2}^{s} \hookrightarrow F_{r,2}^{s} \hookrightarrow W^{s,r} \hookrightarrow \mpq$, so when $r\leq 2$, we have $B_{r,r}^{s} \hookrightarrow \mpq$, by Theorem \ref{besov to mpq}, we know that if $r\leq q$, then $s\geq \tau(r,q)$; if $r>q$, then $s>\tau(r,q)$. 
		
		When $r>2$, we have $B_{r,2}^{s} \hookrightarrow \mpq$, also by  Theorem \ref{besov to mpq}, we know that if $2\leq q,$ then $s\geq \tau(r,q)$; if $2>q$, then $s>\tau(r,q)$.
		
		\item[(iii)] We prove that $W^{\tau(r,q),r} \hookrightarrow  \mpq$ is not true for any $2\leq q <r \leq p$, in which case $\tau(r,q)= d(1/q-1/r)$, for convenience we also denote it by $\tau$.
		
		If not, we have $W^{\tau,r} \hookrightarrow  \mpq$ which is equivalent to $L^{r} \hookrightarrow \mpq^{-\tau}$. So we have \begin{align} \label{eq lr to mpqtau}
			\norm{f}_{\mpq^{-\tau}} \lesssim \norm{f}_{r}.
		\end{align}
		Choose $\eta \in \sch$, such that supp $\hat{\eta} \subset [-1/8,1/8]^{d}$, take $f(x) = \sum_{k\in\Z^{d}} a_{k} e^{ikx} \eta(x-k)$, then $\hat{f}(\xi) = \sum_{k\in\Z^{d}} a_{k} e^{ik(\xi-k)} \hat{\eta}(\xi-k)$. Then by the orthogonality of $\Box_{k}$, we have \begin{align*}
			\Box_{k'}f=\begin{cases*}
				a_{k} e^{ikx} \eta(x-k); & $k=k'$;\\
				0; & $k\neq k'$. 		
			\end{cases*}
		\end{align*}
		So, we have $\norm{f}_{\mpq^{-\tau}} \approx \norm{\inner{k}^{-\tau} a_{k}}_{\ell_{k}^{q}}$.
		
		By the orthogonality of $L^{2}$ space, we know that 
		\begin{align*}
			\norm{f}_{2} = \left(\sum_{k\in\Z^{d}} \norm{a_{k} e^{ikx} \eta(x-k)}_{2}^{2} \right)^{1/2} \approx \norm{a_{k}}_{\ell^{2}}.
		\end{align*}
		By the rapidly decrease of $\eta$, we know that 
		\begin{align*}
			\abs{f(x)} \leq \sum_{k\in\Z^{d}} \abs{a_{k}} \abs{\eta(x-k)} \leq \norm{a_{k}}_{\ell_{k}^{\infty}} \sum_{k\in\Z^{d}} (1+\abs{x-k})^{-N} \lesssim \norm{a_{k}}_{\ell_{k}^{\infty}}.
		\end{align*}
		So we have $\norm{f}_{\infty} \lesssim \norm{a_{k}}_{\ell_{k}^{\infty}}$, then by interpolation, we have $\norm{f}_{r} \lesssim \norm{a_{k}}_{\ell_{k}^{r}}$.
		
		Take the two estimates into \eqref{eq lr to mpqtau}, we have $\norm{\inner{k}^{-\tau} a_{k}} \lesssim \norm{a_{k}}_{\ell_{k}^{r}}$, which is equivalent to \begin{align*}
			\sum_{k\in\Z^{d}} \inner{k}^{-\tau q} b_{k} \lesssim \left( \sum_{k\in\Z^{d}} b_{k}^{r/q}\right)^{q/r}.
		\end{align*}
		So we know that $\set{\inner{k}^{-\tau q}} \in \ell^{(r/q)'}$, which is equivalent to $\tau q (r/q)' >d$, so we have $\tau >d(1/q-1/r)$, which is a contradiction.
		
		\item[(iv)] We prove that when $r=1\leq p, 0< q <\infty,s=\tau(1,q)=d/q$, $W^{s,1} \hookrightarrow M_{p,q}$ is not true.
		
		If not, we have $W^{d/q,1} \hookrightarrow M_{p,q} \hookrightarrow M_{\infty,q}$ which is equivalent to $L^{1} \hookrightarrow M_{\infty,q}^{-d/q}$, i.e. \begin{align} \label{eq L1 to mpq}
			\norm{f}_{M_{\infty,q}^{-d/q}} \lesssim \norm{f}_{1}.
		\end{align}
		Choose $\eta \in \sch$ such that $\hat{\eta} = 1 $ on $[-1,1]^{d}$, for $0<t<1$, take $f(x)=t^{-d} \eta(x/t)$, then $\hat{f}=1 $ on $\frac{1}{t}[-1,1]^{d}$, denote $K_{t}= \set{ k\in \Z^{d}: k+[-3/4,3/4]^{d} \subseteq \frac{1}{t}[-1,1]^{d}}$, then we have \begin{align*}
			\norm{f}_{M_{\infty,q}^{-d/q}} = \norm{\inner{k}^{-d/q} \norm{\Box_{k}f}_{\infty}}_{\ell_{k}^{q}} \geq \norm{\inner{k}^{-d/q} \norm{\Box_{k}f}_{\infty}}_{\ell_{k\in K_{t}}^{q}} \approx \norm{\inner{k}^{-d/q}}_{\ell_{k\in K_{t}}^{q}}.
		\end{align*}
		Also, $\norm{f}_{1} = \norm{\eta}_{1} \approx 1$. So, take the two estimates into \eqref{eq L1 to mpq}, we have $ \norm{\inner{k}^{-d/q}}_{\ell_{k\in K_{t}}^{q}}\lesssim 1$, take $t\rightarrow 0$, we have $\norm{\inner{k}^{-d/q}}_{\ell_{k\in \Z^{d}}^{q}} \lesssim 1$, which is a contraction.
	\end{description}
	
	In fact, we already get the result as desired. To be specific,from (i,ii) we can get the necessity of the condition (2)(3); the necessity of condition (4) can be got by (ii,iv); as for condition (1), when $q<r\leq 2$ or $r>2>q$, the necessity is proved in (ii), when $2\leq q <r$, the necessity is proved in (iii).
	
\subsection{Proof of Theorem \ref{thm mpq to wsr}}

When $1\leq q \leq \infty$, we can get the result as desired just by duality of Theorem \ref{thm wsr to mpq}, so we only need to consider the case $0<q<1$, which could only happen in case (2) or (3). Also, we know that $\sigma(r,q)=0$ in these cases.

Sufficiency: by Lemma \ref{lemma mpq} and dual of Proposition \ref{prop wsp to mpq}, we can get 
$\mpq \hookrightarrow M_{r,1} \hookrightarrow W^{s,r}$ when $s\leq 0$.

Necissity: if we have $\mpq \hookrightarrow W^{s,r}$ which is equivalent to $\mpq^{-s} \hookrightarrow L^{r}$. So, we have \begin{align*}
\norm{f}_{r} \lesssim \norm{f}_{\mpq^{-s}}.
\end{align*}
	For any $k\in \Z^{d}$, take $f(x)=e^{ikx} \eta(x)$ into the equation above, where $\eta \in \sch$ and supp $\hat{\eta} \subseteq[-1/8,1/8]^{d}$. Then we have $1\lesssim \inner{k}^{-s}$, which means that $s\leq 0$.
	
	\section{Embedding between Triebel spaces and modulation spaces}\ 
	
	In this section, we mainly study the embedding $F_{p_{0},q}^{s} \hookrightarrow \mpq$. We first consider a special case like $F_{p,2}^{s} \hookrightarrow \mpq$, which is just the same as Proposition \ref{prop wsp to mpq} in some cases. And using this embedding and interpolation method we can get the result we need.
	
	As for special case, we have:
	\begin{thm}
		\label{thm Fp2 to mpq}
		Let $0< p,q \leq \infty$, then $F_{p,2}^{s} \hookrightarrow \mpq$ if and only if one of the following conditions is satisfied:
		\begin{description}
			\item[(1)] $q\geq p, s \geq \tau(p,q)$;
			\item[(2)] $q<p, s> \tau(p,q)$.
		\end{description}
	\end{thm}
	
	As for the other side of this kind of embedding, we also have:
	\begin{thm}\label{thm mpq to Fp2s}
		Let $0< p,q \leq \infty$, then $\mpq \hookrightarrow F_{p,2}^{s}$ if and only if one of the following conditions is satisfied :
		\begin{description}
			\item[(1)] $q\leq p, s\leq \sigma(p,q)$;
			\item[(2)] $q>p,s< \sigma(p,q)$.
		\end{description}
	\end{thm}

	Then, our main result is:
	\begin{thm}
		\label{thm Fpq to mpq}
		Let $0< p,p_{0},q \leq \infty$, then $F_{p_{0},q}^{s} \hookrightarrow \mpq$ if and only if $p_{0} \leq p$ and one of the following conditions is satisfied:
		\begin{description}
			\item[(1)] $ p_{0} \leq q, s\geq \tau(p_{0},q)$;
			\item[(2)] $ p_{0} >q, s> \tau(p_{0},q)$.
		\end{description}
	\end{thm}
	
	\begin{figure}
		\centering
		\begin{tikzpicture}[scale=5]
			\filldraw[fill=gray!20] (0,0) --(1.5,1.5) --(0,1.5);
			\filldraw[fill=gray!40] (0,0)--(1.5,1.5)--(1.5,0);
		%	\filldraw[fill = gray!60] (0.5,0.5)--(1,0)--(1.5,0)--(1.5,1.5)--(0.5,1.5);
			\draw (0,0)--(1.5,1.5);
			\draw [->] (0,0) -- (1.55,0);
			\draw[->] (0,0) -- (0,1.55);
			\draw [dashed] (0.5,0.5)--(1,0);
			\draw [dashed] (0.5,0.5)--(0.5,1.5);
			\draw[- >] (0.6,0.6) -- (0.7,0.5) ;
			\draw (1/3,2/3) node {(1)} ;
			\draw (0.8,0.4) node {(2)};
			\draw (1.55,-2pt) node {$\frac{1}{p_{0}}$};
			\draw (-2pt,1.55) node {$\frac{1}{q}$};
			\draw (-1pt,-1pt) node {0};
			\draw (1,-2pt) node {1};
			\draw (0.25,1.25) node {$\tau= d(\frac{1}{q}-\frac{1}{p_{0}})$};
			\draw (0.5,0.2) node {$\tau =0$};
			\draw (1.1,1.1) node {$\tau = d(\frac{1}{p_{0}} +\frac{1}{q}-1)$};

		\end{tikzpicture}
		\caption{the two cases of $(1/p_{0},1/q)$ in Theorem \ref{thm Fpq to mpq}}
	\end{figure}
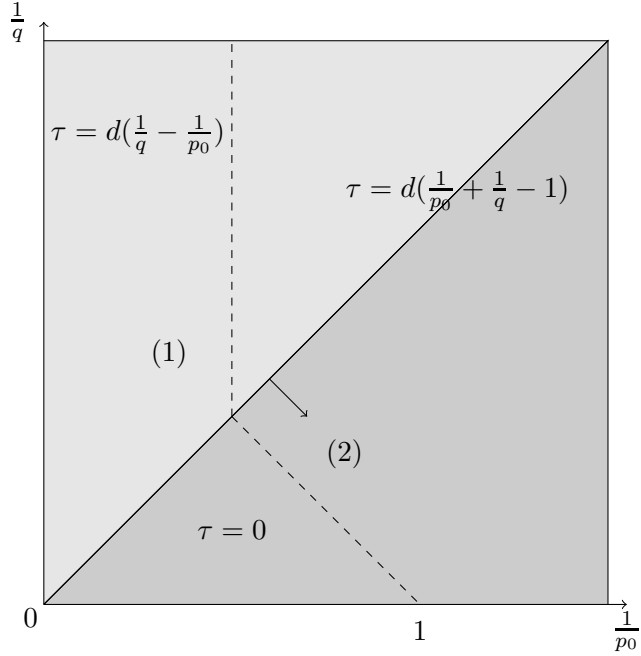

	As for the other side of this kind of embedding, we also have:
	\begin{thm} \label{thm mpq to Frq}
		Let $0< p,p_{1},q \leq \infty$, then $\mpq \hookrightarrow F_{p_{1},q}^{s}$ if and only if $p_{1} \geq p$ and one of the following conditions is satisfied:
		\begin{description}
			\item[(1)] $p_{1} \geq q, s\leq \sigma(p_{1},q)$;
			\item[(2)] $p_{1} <q,s< \sigma(p_{1},q)$.
		\end{description}
	\end{thm}

	\subsection{Proof of Theorem \ref{thm Fp2 to mpq}}\ 
	
	For $1<p<\infty$, we know that $F_{p,2}^{s} \approx W^{s,p}$, so by Theorem \ref{thm wsr to mpq}, we get the result as desired.
	
	For $0<p\leq 1$, we know that $F_{p,2}^{s}=h_{p}$ the local Hardy space, so by the result in \cite{Kobayashi2009Embedding}, we get the result as desired.
	
	For $p=\infty$, if we have $F_{\infty,2}^{s} \hookrightarrow \mpq$, then by Lemma \ref{lemma besov triebel}, we have $W^{\infty,s} \hookrightarrow F_{\infty,2}^{s} \hookrightarrow \mpq$, then by Theorem \ref{thm wsr to mpq}, we have the condition as desired. As for sufficiency, in condition (1), by Lemma \ref{lemma besov triebel}, we have $F_{\infty,2}^{s} \hookrightarrow B_{\infty,\infty}^{0}$, then by Proposition \ref{prop besov to mpq}, we have $B_{\infty,\infty}^{0} \hookrightarrow M_{\infty,\infty}$; in condition (2), by Lemma \ref{lemma besov triebel}, we have $F_{\infty,2}^{s} \hookrightarrow B_{\infty,\infty}^{s} \hookrightarrow B_{\infty,q}^{d/q}$, then by Proposition \ref{prop besov to mpq}, we have $ B_{\infty,q}^{d/q} \hookrightarrow M_{\infty,q}$.
	
	\subsection{Proof of Theorem \ref{thm mpq to Fp2s}}
	By the same discussion in the proof of Theorem \ref{thm Fp2 to mpq}, we only need to consider the case when $p=\infty$. Also, when $1\leq q \leq \infty$,  by duality, we can get the result as desired. When $0<q<1,p=\infty$, we have $\sigma(p,q) =0$. 
	
	Sufficiency: by Lemma \ref{lemma mpq}, we have $M_{\infty,q} \hookrightarrow M_{\infty,1} $. By Theorem \ref{thm Fp2 to mpq}, we have $F_{1,2} \hookrightarrow M_{1,\infty}$, by duality, we have $M_{\infty,1} \hookrightarrow F_{\infty,2}$, so we have $M_{\infty,q} \hookrightarrow F_{\infty,2} \hookrightarrow F_{\infty,2}^{s}$ for $s\leq 0$.
	
	Necessity: if we have $M_{\infty,q} \hookrightarrow F_{\infty,2}^{s}$, then by Lemma \ref{lemma besov triebel}, we have $F_{\infty,2}^{s} \hookrightarrow B_{\infty,\infty}^{s}$. By Theorem \ref{mpq to besov}, we have $s\leq \sigma(\infty,q)$.

	\subsection{Proof of Theorem \ref{thm Fpq to mpq}}
	
	\subsubsection{Sufficiency}\ 
	
	For condition (1): by Lemma \ref{lemma besov triebel}, we have $F_{p_{0},q}^{s} \hookrightarrow B_{p_{0},q}^{s}$, then by Theorem \ref{besov to mpq}, we have $B_{p_{0},q}^{s} \hookrightarrow \mpq$.
	
	For condition (2): by Lemma \ref{lemma besov triebel}, we have $F_{p_{0},q}^{s} \hookrightarrow B_{p_{0},p_{0}}^{s} \hookrightarrow B_{p_{0},q}^{\tau(p_{0},q)+\eps}$ for some $0<\eps <<1$, then by Theorem \ref{besov to mpq}, we have $B_{p_{0},q}^{\tau(p_{0},q)+\eps} \hookrightarrow \mpq$.
	
	\subsubsection{Necessity}
	If  we have $F_{p_{0},q}^{s} \hookrightarrow \mpq$, then we have \begin{align}
		\label{eq fp0q to mpq}
		\norm{f}_{\mpq} \lesssim \norm{f}_{F_{p_{0},q}^{s}}.
	\end{align}
	\begin{description}
		\item[(i)] By the same method as in subsection \ref{subsubsec Wrq to mpq}, we can get $p_{0} \leq p$.
		\item[(ii)] By Lemma \ref{lemma besov triebel}, we have $B_{p_{0},p_{0}\wedge q}^{s} \hookrightarrow F_{p_{0},q}^{s} \hookrightarrow \mpq$, then by Theorem \ref{besov to mpq}, we have $s\geq \tau(p_{0},q)$.
		\item[(iii)] When $p_{0} >q \geq 2$, by Lemma \ref{lemma besov triebel}, we have $W^{s,p_{0}} \approx F_{p_{0},2}^{s} \hookrightarrow F_{p_{0},q}^{s} \hookrightarrow \mpq$, then by Theorem \ref{thm wsr to mpq}, we have $s>\tau(p_{0},q)$.
		\item[(iv)] When $2\geq p_{0} >q$, by a similar proof of Proposition 3.2 in \cite{Guo2017Inclusion}, we have $\norm{2^{jd/q} a_{j}}_{\ell_{j}^{q}} \lesssim \norm{2^{j(s+d(1-1/p_{0}))} a_{j}}_{p_{0}}$, then we have $s>d(1/q+1/p_{0}-1)$.
		\item[(v)] We prove that $F_{p_{0},q}^{\tau(p_{0},q)} \hookrightarrow \mpq$ is not true for any $p_{0}>q, q<2, p_{0}>2$, in which case $\tau(p_{0},q) =d(1/q-1/p_{0})$. 
		
		If not, we have $F_{p_{0},q}^{d(1/q-1/p_{0})} \hookrightarrow \mpq$, also by sufficiency, we have $F_{p_{0},p_{0}}^0\hookrightarrow M_{p_{0},p_{0}} \hookrightarrow M_{p,p_{0}}$, then by interpolation in Lemma \ref{lemma besov triebel} and \ref{lemma mpq}, we have $F_{p_{0},2}^{d(1/2-1/p_{0})} \hookrightarrow M_{p,2}$ which is contradiction with case (iii). 
	\end{description}
	
	Combine the case(i),(ii), we get the condition (1); combine the case (i),(iii),(iv),(v), we get the condition (2) as desired.
	
	\begin{rem}\ 
	\begin{enumerate}
	\item As for case (v) in the proof above, we prove it by contradiction using the method of interpolation. Also, we can prove it directly by take some good function into the norm inequality. One can modify the proof of Proposition 3.1 in \cite{Guo2017Inclusion} to get the result.
	\item The proof of Theorem \ref{thm mpq to Frq} is similar to the proof above. For the reader's convenience, we give the outline of the proof here. The sufficiency can get by the relation of Triebel spaces and Besov spaces, then use Theorem \ref{mpq to besov} can get the result. As for necessity, we can use a revised version of Proposition 3.1,3.2 in \cite{Guo2017Inclusion} directly to get the result as desired.
	\item  As we can see, the embedding we consider above is not the most general case like $F_{p_{0},q_{0}}^{s} \hookrightarrow \mpq$. I believe that for the general case the sufficient and necessary is the similar to Theorem \ref{thm Fpq to mpq}. But I have not proven it yet. The most difficult case left is the same as the open problem in \cite{Guo2017Inclusion}.
	\end{enumerate}
	\end{rem}

	\section{Embedding between Fourier $L^{p}$ spaces and modulation spaces}\ 
	
	In this section, we consider the embedding $\mpq^{s} \hookrightarrow \F L^{r}$. Because of the orthogonality of $\Box_{k}(\sigma_{k})$, we can calculate the $\F L^{r}$ norm easily. Then, we can get the result as follow. Recall that when $0<p<1$, we denote $p'=\infty$.

	\begin{thm}
		\label{thm mpq to Flp}
		Let $0< p,q,r\leq \infty,s\in \real$, then we have $\mpq^{s} \hookrightarrow \F L^{r}$ if and only if one of the following conditions is satisfied:
		\begin{description}
			\item[(1)] $p\leq 2, q\leq r \leq p',s\geq 0$;
			\item[(2)] $p\leq 2, r\leq p', r<q, s>d(1/r-1/q)$.
		\end{description}	
	\end{thm}
	
	As for the other side of this kind of embedding, we also have:
	
	\begin{thm}
		\label{thm Flp to mpq}
		Let Let $0< p,q,r\leq \infty,s\in \real$, then we have $ \F L^{r} \hookrightarrow \mpq^{s}$ if and only if one of the following conditions is satisfied:
		\begin{description}
			\item[(1)] $p\geq2,p'\leq r \leq q, s\leq 0$;
			\item[(2)] $ p\geq 2, p'\leq r,r>q, s<d(1/r-1/q)$.
		\end{description}
	\end{thm}
	
	\subsection{Proof of Theorem \ref{thm mpq to Flp}}
	
	\subsubsection{Sufficiency}\ 
	
	If we have condition (1), when $1\leq p \leq \infty$, by the orthogonality of $\sigma_{k}$, H$\ddot{o}$lder's inequality of $L^p$ spaces and Hausdorff-Young's inequality of Fourier transform, we have \begin{align*}
		\norm{\hat{f}}_{r} \approx \norm{ \norm{\sigma_{k} \hat{f}}_{r}}_{\ell_{k}^{r}} \leq \norm{ \norm{\sigma_{k} \hat{f}}_{p'}}_{\ell_{k}^{r}} \leq \norm{\norm{\Box_{k}f}_{p}}_{\ell_{k}^{r}} \leq \norm{\norm{\Box_{k}f}_{p}}_{\ell_{k}^{q}} = \norm{f}_{\mpq} \leq \norm{f}_{\mpq^{s}}.
	\end{align*}
	
	When $0<p<1$, by Lemma \ref{lemm bernstein}, we have $\norm{\sigma_{k} \hat{f}}_{p'} \leq \norm{\Box_{k}f}_{1} \leq C \norm{\Box_{k}f}_{p}$, then by the same estimates we can get the result as desired.
	
	If we have condition (2), by Lemma \ref{lemma mpq}, we have $\mpq^{s} \hookrightarrow M_{p,r}$, then by the sufficiency of condition(1), we have $M_{p,r} \hookrightarrow \F L^{r}$.

	\subsubsection{Necessity}
	
	If we have $\mpq^{s} \hookrightarrow \F L^{r}$, then we have \begin{align}
		\label{eq mpq to Flr}
		\norm{\hat{f}}_{r} \lesssim \norm{f}_{\mpq^{s}}.
	\end{align}
	The conditions we need are contained the following cases.
	\begin{description}
		\item[(i)] For any $f \in \sch$ with supp $\hat{f} \subseteq [-1/8,1/8]^{d}$, in this case we have $\norm{f}_{\mpq^{s}} \approx \norm{f}_{p}$. So, \eqref{eq mpq to Flr} means that \begin{align*}
			\norm{\hat{f}}_{r} \lesssim \norm{f}_{p}.
		\end{align*} 
		Then by Hausdorff-Young's inequality on compact set (see Tao's note of Math 254B, one can extend the result to the case of $0<p<1$ with the same proof), we have $r\leq p',p\leq 2$.
		\item[(ii)] For any $k\in \Z^{d}$, take $f(x)=e^{ikx} \eta(x)$, where $\eta \in \sch,$supp  $\hat{\eta} \subseteq [-1/8,1/8]^{d}$, then we know that supp $\hat{f} \subseteq k+[-1/8,1/8]^{d}$. So, we have \begin{align*}
			\norm{f}_{\mpq^{s}} = \inner{k}^{s} \norm{\eta}_{p} \approx \inner{k}^{s}; \ \ \norm{f}_{\F L^{r}} = \norm{\hat{\eta}}_{r} \approx 1.
		\end{align*}
		Take the estimates into \eqref{eq mpq to Flr}, we have $1\lesssim \inner{k}^{s}$, which means that $s\geq 0$, which is the result in condition (1).
		\item[(iii)] When $r<q$, take $f(x) = \sum_{k\in\Z^{d}} a_{k} e^{ikx} \eta(x)$, $\eta$ is the same as case (ii). Then we have $\hat{f} (\xi) = \sum_{k\in\Z^{d}} a_{k} \hat{\eta}(\xi-k)$. 
		So, we have \begin{align*}
			\norm{f}_{\mpq^{s}} \approx \norm{\inner{k}^{s} a_{k}}_{\ell_{k}^{q}}; \ \ \norm{\hat{f}}_{r} \approx \norm{a_{k}} _{\ell_{k}^{r}}.
		\end{align*}
		Take the estimates above into \eqref{eq mpq to Flr}, we have $\norm{a_{k}}_{\ell_{k}^{r}} \lesssim \norm{\inner{k}^{s} a_{k}}_{\ell_{k}^{q}}$, then we have $\norm{\inner{k}^{-sr} b_{k}}_{\ell_{k}^{1}} \lesssim \norm{b_{k}}_{\ell_{k}^{q/r}}$, which means that $\set{\inner{k}^{-sr}} \in \ell^{(q/r)'}$. So, we have $sr(q/r)' >d$, i.e. $s>d(1/r-1/q)$, which is the result in condition (2).
	\end{description}
	
	\begin{rem}
	The proof of Theorem \ref{thm Flp to mpq} is similar to the proof above, we omit it here.
	\end{rem}

%	\section*{References}
	\bibliographystyle{ieeetr} 	
	\bibliography{besovmpq2}

\end{document}